\documentclass[a4paper,11pt,reqno]{amsart}

\usepackage[utf8]{inputenc}
\usepackage[T1]{fontenc}
\usepackage{latexsym,amsmath,amsfonts,amscd,amssymb,amsthm,mathrsfs}
\usepackage{enumerate}
\usepackage[all]{xy}
\usepackage{color}
\usepackage{tikz}
\usepackage{graphicx}
\usepackage{caption}
\usepackage{subcaption}

\definecolor{cobalt}{RGB}{61,89,171}
\usepackage[colorlinks,citecolor=cobalt,linkcolor=cobalt,urlcolor=cobalt,pdfpagemode=UseNone,backref = page]{hyperref}

\newcommand{\Z}{\mathbb{Z}}
\newcommand{\R}{\mathbb{R}}
\newcommand{\C}{\mathbb{C}}

\newcommand{\tens}[1]{
  \mathbin{\mathop{\otimes}\displaylimits_{#1}}
}

\DeclareMathOperator{\coker}{coker}
\DeclareMathOperator{\im}{im}

\theoremstyle{plain}
\newtheorem{theorem}{Theorem}[section]

\newtheorem{problem}[theorem]{Problem}
\theoremstyle{definition}

\theoremstyle{remark}
\newtheorem{remark}[theorem]{Remark}

\begin{document}

\title{Hodge numbers of a hypothetical complex structure on $S^6$}

\author[D. Angella]{Daniele Angella}

\address[D. Angella]{
Dipartimento di Matematica e Informatica "Ulisse Dini"\\
Università di Firenze\\
viale Morgagni 67/a\\
50134 Firenze\\
Italy\\
}

\email{daniele.angella@gmail.com}
\email{daniele.angella@unifi.it}

\urladdr{http://sites.google.com/site/danieleangella/}

\thanks{The author is supported by the Project SIR2014 AnHyC ``Analytic aspects in complex and hypercomplex geometry'' (code RBSI14DYEB), and by GNSAGA of INdAM}

\date{\today}

\begin{abstract}
These are the notes for the talk "Hodge numbers of a hypothetical complex structure on $S^6$" given by the author at the MAM1 ``(Non)-existence of complex structures on $S^6$'' held in Marburg in March 2017. They are based on \cite{gray} and \cite{ugarte}, where Hodge numbers and the dimensions of the successive pages of the Fr\"olicher spectral sequence for $S^6$ endowed with a hypothetical complex structure are investigated. We also add results from \cite{mchugh-1}, where the Bott-Chern cohomology of hypothetical complex structures on $S^6$ is studied.
The material is not intended to be original.
\end{abstract}

\maketitle

\section*{Introduction}

In \cite{gray}, Alfred Gray investigated the properties of the Dolbeault cohomology of a hypothetical complex structure on $S^6$. In particular, he proved that, if $S^6$ admits a complex structure, then the corresponding Hodge number $h^{0,1}$ is positive \cite[Theorem 5]{gray}, that is, there exists a non-zero $\overline\partial$-closed non-$\overline\partial$-exact $(0,1)$-form, that hopefully may be interpreted geometrically. In particular, the Fr\"olicher spectral sequence does not degenerate at the first page. Moreover, one can deduce also that, at the second degree, (that is, for $h^{p,q}$ with $p+q=2$,) either $h^{2,0}$ or $h^{1,1}$ are non-zero \cite[Proposition 3.1]{ugarte}. As in \cite[Proposition 10.3]{huckleberry-kebekus-peternell}, \cite[Lemma 1, Lemma 2]{mchugh-1}, one can also show that $h^{1,0}(X)\leq h^{2,0}(X)$ (and that $h^{1,0}(X)\leq1$ once proven that the algebraic dimension is zero, see \cite{rollenske-mam1}); as in \cite[Remark 3.4]{ugarte}, one can also show that $h^{1,1}(X) \geq h^{1,2}(X) - h^{0,2}(X)$.

Further invariants, besides Dolbeault cohomology, are provided by the successive pages of the Fr\"olicher spectral sequence. This has been investigated by Luis Ugarte in \cite{ugarte}. In particular, he proves that the Fr\"olicher spectral sequence degenerates at the third page \cite[Lemma 2.1]{ugarte}, that Serre symmetry holds also for the successive pages of the Fr\"olicher spectral sequence \cite[page 174]{ugarte}, and that either $h^{1,1}$ is non-zero, or there exists a $d$-closed holomorphic $2$-form and $E_2 \neq E_3 = E_\infty$ \cite[Corollary 3.3.]{ugarte}; this latter special case is depicted in Figure \ref{fig:frolicher-h11zero}.
More precisely, the Fr\"olicher spectral sequence degenerates at the second page if and only if $h^{1,2}(X)=h^{1,1}(X)-1\geq 0$ and $h_2^{0,2}(X)=0$.
Note that this would be the first simply-connected example in the lowest possible dimension with $E_2 \neq E_\infty$, as asked in \cite[page 444]{griffiths-harris}, see \cite{pittie, cordero-fernandez-ugarte-gray, rollenske, bigalke-rollenske} (for non-simply-connected example, see \cite{cordero-fernandez-gray}). 

\begin{figure}
    \centering
    \begin{subfigure}[t]{0.25\textwidth}

\begin{tikzpicture}
\newcommand\un{.8}

\draw[help lines, step=\un] (0,0) grid (4*\un,4*\un);

\begingroup\makeatletter\def\f@size{8}\check@mathfonts
\node at (2*\un,4.3*\un) {$H^{\bullet,\bullet}_{\overline\partial}(X)$};
\endgroup

\begingroup\makeatletter\def\f@size{7}\check@mathfonts
\foreach \x in {0,...,3}
  \node at (\un*.5+\un*\x,-.3) {\x};
\foreach \y in {0,...,3}
  \node at (-.3,\un*.5+\un*\y) {\y};
\endgroup

\coordinate (A00) at (0*\un+1/2*\un, 0*\un+1/2*\un);
\coordinate (A01) at (0*\un+1/2*\un, 1*\un+1/2*\un);
\coordinate (A10) at (1*\un+1/2*\un, 0*\un+1/2*\un);
\coordinate (A20) at (2*\un+1/2*\un, 0*\un+1/2*\un);
\coordinate (A11) at (1*\un+1/2*\un, 1*\un+1/2*\un);
\coordinate (A02) at (0*\un+1/2*\un, 2*\un+1/2*\un);
\coordinate (A30) at (3*\un+1/2*\un, 0*\un+1/2*\un);
\coordinate (A21) at (2*\un+1/2*\un, 1*\un+1/2*\un);
\coordinate (A12) at (1*\un+1/2*\un, 2*\un+1/2*\un);
\coordinate (A03) at (0*\un+1/2*\un, 3*\un+1/2*\un);
\coordinate (A31) at (3*\un+1/2*\un, 1*\un+1/2*\un);
\coordinate (A22) at (2*\un+1/2*\un, 2*\un+1/2*\un);
\coordinate (A13) at (1*\un+1/2*\un, 3*\un+1/2*\un);
\coordinate (A32) at (3*\un+1/2*\un, 2*\un+1/2*\un);
\coordinate (A23) at (2*\un+1/2*\un, 3*\un+1/2*\un);
\coordinate (A33) at (3*\un+1/2*\un, 3*\un+1/2*\un);

\begingroup\makeatletter\def\f@size{4}\check@mathfonts
\node at (A00) {$1$};
\node at (A01) {$h^{0,2}+1$};
\node at (A10) {$h^{1,0}$};
\node at (A20) {$h^{2,0}$};
\node at (A11) {$h^{1,1}$};
\node at (A02) {$h^{0,2}$};
\node at (A30) {$0$};
\node at (A21) {$h^{1,2}$};
\node at (A12) {$h^{1,2}$};
\node at (A03) {$0$};
\node at (A31) {$h^{0,2}$};
\node at (A22) {$h^{1,1}$};
\node at (A13) {$h^{2,0}$};
\node at (A32) {$h^{0,2}+1$};
\node at (A23) {$h^{1,0}$};
\node at (A33) {$1$};
\endgroup

\end{tikzpicture}
\caption{First page, namely, Dolbeault cohomology.}
\label{fig:dolbeault}
    \end{subfigure}
    \qquad
    \begin{subfigure}[t]{0.25\textwidth}
\begin{tikzpicture}
\newcommand\un{.8}

\draw[help lines, step=\un] (0,0) grid (4*\un,4*\un);

\begingroup\makeatletter\def\f@size{8}\check@mathfonts
\node at (2*\un,4.3*\un) {$E^{\bullet,\bullet}_{2}(X)$};
\endgroup

\begingroup\makeatletter\def\f@size{7}\check@mathfonts
\foreach \x in {0,...,3}
  \node at (\un*.5+\un*\x,-.3) {\x};
\foreach \y in {0,...,3}
  \node at (-.3,\un*.5+\un*\y) {\y};
\endgroup

\coordinate (A00) at (0*\un+1/2*\un, 0*\un+1/2*\un);
\coordinate (A01) at (0*\un+1/2*\un, 1*\un+1/2*\un);
\coordinate (A10) at (1*\un+1/2*\un, 0*\un+1/2*\un);
\coordinate (A20) at (2*\un+1/2*\un, 0*\un+1/2*\un);
\coordinate (A11) at (1*\un+1/2*\un, 1*\un+1/2*\un);
\coordinate (A02) at (0*\un+1/2*\un, 2*\un+1/2*\un);
\coordinate (A30) at (3*\un+1/2*\un, 0*\un+1/2*\un);
\coordinate (A21) at (2*\un+1/2*\un, 1*\un+1/2*\un);
\coordinate (A12) at (1*\un+1/2*\un, 2*\un+1/2*\un);
\coordinate (A03) at (0*\un+1/2*\un, 3*\un+1/2*\un);
\coordinate (A31) at (3*\un+1/2*\un, 1*\un+1/2*\un);
\coordinate (A22) at (2*\un+1/2*\un, 2*\un+1/2*\un);
\coordinate (A13) at (1*\un+1/2*\un, 3*\un+1/2*\un);
\coordinate (A32) at (3*\un+1/2*\un, 2*\un+1/2*\un);
\coordinate (A23) at (2*\un+1/2*\un, 3*\un+1/2*\un);
\coordinate (A33) at (3*\un+1/2*\un, 3*\un+1/2*\un);

\begingroup\makeatletter\def\f@size{4}\check@mathfonts
\node at (A00) {$1$};
\node at (A01) {$h^{0,1}_2$};
\node at (A10) {$0$};
\node at (A20) {$h^{0,1}_2$};
\node at (A11) {$0$};
\node at (A02) {$h^{0,2}_2$};
\node at (A30) {$0$};
\node at (A21) {$h^{0,2}_2$};
\node at (A12) {$h^{0,2}_2$};
\node at (A03) {$0$};
\node at (A31) {$h^{0,2}_2$};
\node at (A22) {$0$};
\node at (A13) {$h^{0,1}_2$};
\node at (A32) {$h^{0,1}_2$};
\node at (A23) {$0$};
\node at (A33) {$1$};
\endgroup

\end{tikzpicture}
\caption{Second page.}
\label{fig:E2}
    \end{subfigure}
    \qquad
    \begin{subfigure}[t]{0.25\textwidth}
\begin{tikzpicture}
\newcommand\un{.8}

\draw[help lines, step=\un] (0,0) grid (4*\un,4*\un);

\begingroup\makeatletter\def\f@size{8}\check@mathfonts
\node at (2*\un,4.7*\un) {$E^{\bullet,\bullet}_{r}(X)$};
\endgroup
\begingroup\makeatletter\def\f@size{6}\check@mathfonts
\node at (2*\un,4.3*\un) {with $r\geq 3$};
\endgroup

\begingroup\makeatletter\def\f@size{4}\check@mathfonts
\foreach \x in {0,...,3}
  \node at (\un*.5+\un*\x,-.3) {\x};
\foreach \y in {0,...,3}
  \node at (-.3,\un*.5+\un*\y) {\y};
\endgroup

\coordinate (A00) at (0*\un+1/2*\un, 0*\un+1/2*\un);
\coordinate (A01) at (0*\un+1/2*\un, 1*\un+1/2*\un);
\coordinate (A10) at (1*\un+1/2*\un, 0*\un+1/2*\un);
\coordinate (A20) at (2*\un+1/2*\un, 0*\un+1/2*\un);
\coordinate (A11) at (1*\un+1/2*\un, 1*\un+1/2*\un);
\coordinate (A02) at (0*\un+1/2*\un, 2*\un+1/2*\un);
\coordinate (A30) at (3*\un+1/2*\un, 0*\un+1/2*\un);
\coordinate (A21) at (2*\un+1/2*\un, 1*\un+1/2*\un);
\coordinate (A12) at (1*\un+1/2*\un, 2*\un+1/2*\un);
\coordinate (A03) at (0*\un+1/2*\un, 3*\un+1/2*\un);
\coordinate (A31) at (3*\un+1/2*\un, 1*\un+1/2*\un);
\coordinate (A22) at (2*\un+1/2*\un, 2*\un+1/2*\un);
\coordinate (A13) at (1*\un+1/2*\un, 3*\un+1/2*\un);
\coordinate (A32) at (3*\un+1/2*\un, 2*\un+1/2*\un);
\coordinate (A23) at (2*\un+1/2*\un, 3*\un+1/2*\un);
\coordinate (A33) at (3*\un+1/2*\un, 3*\un+1/2*\un);

\begingroup\makeatletter\def\f@size{7}\check@mathfonts
\node at (A00) {$1$};
\node at (A01) {$0$};
\node at (A10) {$0$};
\node at (A20) {$0$};
\node at (A11) {$0$};
\node at (A02) {$0$};
\node at (A30) {$0$};
\node at (A21) {$0$};
\node at (A12) {$0$};
\node at (A03) {$0$};
\node at (A31) {$0$};
\node at (A22) {$0$};
\node at (A13) {$0$};
\node at (A32) {$0$};
\node at (A23) {$0$};
\node at (A33) {$1$};
\endgroup

\end{tikzpicture}
\caption{Successive pages.}
\label{fig:E3}
    \end{subfigure}
\caption{Dimensions of the Fr\"olicher spectral sequence of a hypothetical complex structure on $S^6$. With respect to the bi-graduation $\wedge^{p,q}X$, the horizontal axis represents $p$ and the vertical axis represents $q$.
The $h^{0,2}$, $h^{1,0}$, $h^{1,1}$, $h^{1,2}$, $h^{2,0}$, $h^{0,1}_2$, $h^{0,2}_2$ are unknown non-negative integers satisfying $h^{2,0}+h^{1,1}=h^{1,0}+h^{1,2}+1$, $h^{1,0}\leq h^{2,0}$, $h^{1,1} \geq h^{1,2} - h^{0,2}$, and $h_2^{0,1}\leq h^{0,1}$, $h_2^{0,2}\leq h^{0,2}$, $h^{0,1}_2=h^{1,2}-h^{1,1}+1$. (For $h^{1,0}\leq1$, compare\cite{rollenske-mam1}.)}
    \label{fig:frolicher}
\end{figure}
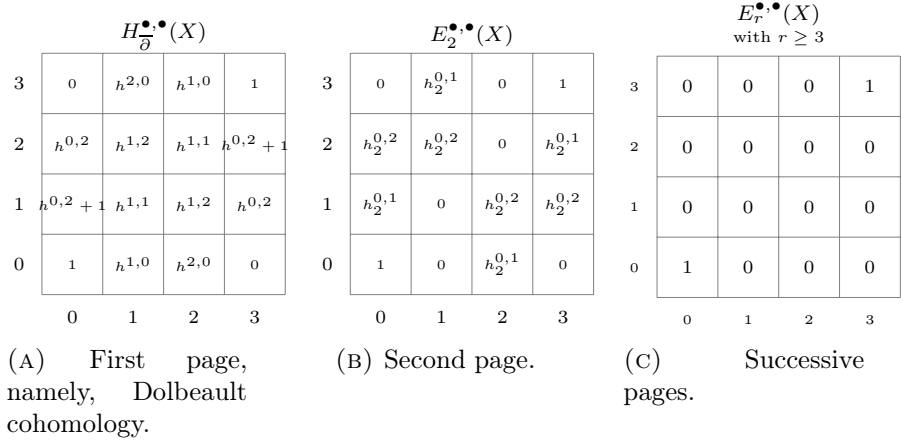

We summarize the informations concerning Dolbeault cohomology \eqref{eq:h00}, \eqref{eq:h30}, \eqref{eq:h01-h02}, \eqref{eq:h20-h11-h10-h12}, respectively the successive pages of the Fr\"olicher spectral sequence \eqref{eq:h203-h230}, \eqref{eq:h200-h233}, \eqref{eq:h210-h223}, \eqref{eq:h222-h211}, \eqref{eq:h2ug2}, \eqref{eq:h2ug3}, \eqref{eq:h2var} in the diamonds in Figure \ref{fig:frolicher}, see also \cite{mchugh-1}.

\subsection{Further remarks}
We add some remarks on the literature.

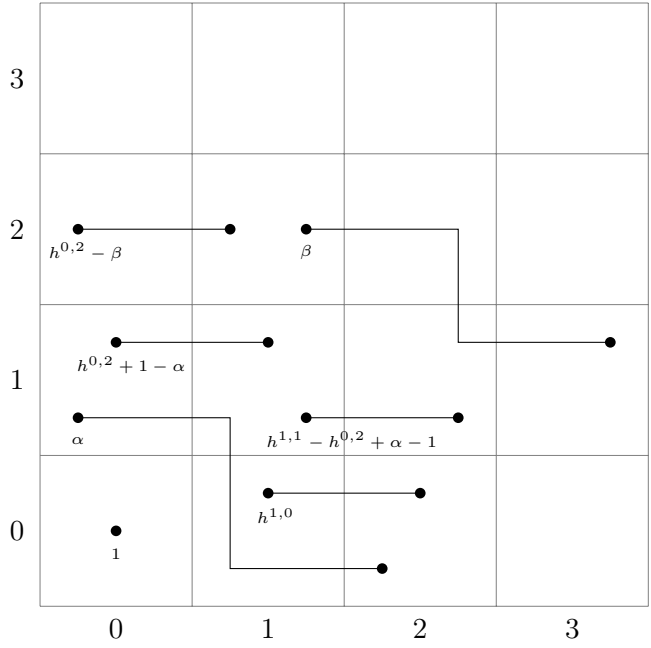
\begin{figure}
\begin{center}
\begin{tikzpicture}
\newcommand\un{2}

\draw[help lines, step=\un] (0,0) grid (4*\un,4*\un);

\foreach \x in {0,...,3}
  \node at (\un*.5+\un*\x,-.3) {\x};
\foreach \y in {0,...,3}
  \node at (-.3,\un*.5+\un*\y) {\y};

\coordinate (A) at (0*\un+1/2*\un, 0*\un+1/2*\un);
\coordinate (B) at (0*\un+1/4*\un, 1*\un+1/4*\un);
\coordinate (C) at (0*\un+1/2*\un, 1*\un+3/4*\un);
\coordinate (D) at (0*\un+1/4*\un, 2*\un+1/2*\un);
\coordinate (E) at (1*\un+1/2*\un, 0*\un+3/4*\un);
\coordinate (F) at (1*\un+1/4*\un, 0*\un+1/4*\un);
\coordinate (G) at (1*\un+1/4*\un, 1*\un+1/4*\un);
\coordinate (H) at (1*\un+3/4*\un, 1*\un+1/4*\un);
\coordinate (I) at (1*\un+1/2*\un, 1*\un+3/4*\un);
\coordinate (L) at (1*\un+1/4*\un, 2*\un+1/2*\un);
\coordinate (M) at (1*\un+3/4*\un, 2*\un+1/2*\un);
\coordinate (N) at (2*\un+1/2*\un, 0*\un+3/4*\un);
\coordinate (O) at (2*\un+1/4*\un, 0*\un+1/4*\un);
\coordinate (P) at (2*\un+3/4*\un, 1*\un+1/4*\un);
\coordinate (Q) at (2*\un+3/4*\un, 1*\un+3/4*\un);
\coordinate (R) at (2*\un+3/4*\un, 2*\un+1/2*\un);
\coordinate (S) at (3*\un+3/4*\un, 1*\un+3/4*\un);

\newcommand{\raggio}{1*\un pt}
\fill (A) circle (\raggio);
\fill (B) circle (\raggio);
\fill (C) circle (\raggio);
\fill (D) circle (\raggio);
\fill (E) circle (\raggio);
\fill (H) circle (\raggio);
\fill (I) circle (\raggio);
\fill (L) circle (\raggio);
\fill (M) circle (\raggio);
\fill (N) circle (\raggio);
\fill (O) circle (\raggio);
\fill (P) circle (\raggio);
\fill (S) circle (\raggio);

\draw (B) -- (G) -- (F) -- (O);
\draw (C) -- (I);
\draw (D) -- (L);
\draw (E) -- (N);
\draw (H) -- (P);
\draw (M) -- (R) -- (Q) -- (S);

\begingroup\makeatletter\def\f@size{6}\check@mathfonts
\node at (0*\un+1/2*\un, 0*\un+1/2*\un-.3) {$1$};
\node at (0*\un+1/4*\un, 1*\un+1/4*\un-.3) {$\alpha$};
\node at (0*\un+1/2*\un+.2, 1*\un+3/4*\un-.3) {$h^{0,2}+1-\alpha$};
\node at (0*\un+1/4*\un+.1, 2*\un+1/2*\un-.3) {$h^{0,2}-\beta$};
\node at (1*\un+1/2*\un+.1, 0*\un+3/4*\un-.3) {$h^{1,0}$};
\node at (1*\un+3/4*\un+.6, 1*\un+1/4*\un-.3) {$h^{1,1}-h^{0,2}+\alpha-1$};
\node at (1*\un+3/4*\un, 2*\un+1/2*\un-.3) {$\beta$};
\endgroup

\end{tikzpicture}
\end{center}
\caption{Main structure of the double complex associated to a hypothetical complex structure on the $6$-dimensional sphere.
The labels count the number of respective objects and $\alpha$, $h^{0,2}$, $\beta$, $h^{1,0}$, $h^{1,1}$ are unknown non-negative integers.
In this diagram, the infinite squares and the arrows arising from symmetries have been removed.}
\label{fig:S6}
\end{figure}

We first notice that, since the de Rham cohomology of $S^6$ is trivial in degree different that $0$ and $6$, from the above informations we can derive the structure of the double-complex of hypothetical complex structures on $S^6$, see Figure \ref{fig:S6} \cite{angella-3}.
Then one gets also informations on the Bott-Chern $H^{\bullet,\bullet}_{BC}(X) := \frac{\ker\partial\cap\ker\overline\partial}{\im\partial\overline\partial}$ and Aeppli $H^{\bullet,\bullet}_{A}(X) := \frac{\ker\partial\overline\partial}{\im\partial+\im\overline\partial}$ cohomologies: they have been investigated by Andrew McHugh in \cite{mchugh-1, mchugh-2}.
More precisely, he obtains that \cite[Section 3]{mchugh-1}, \cite[Proposition 3.5]{mchugh-2}:
$$ h^{0,0}_{BC}(X) = 1 , $$
$$ h^{1,0}_{BC}(X) = h^{0,1}_{BC}(X) = 0 , $$
$$ h^{2,0}_{BC}(X) = h^{0,2}_{BC}(X) = h^{2,0}(X) , $$
$$ h^{1,1}_{BC}(X) = 2\, h^{0,1}(X) , $$
$$ h^{3,0}_{BC}(X) = h^{0,3}_{BC}(X) = 0 ,  $$
$$ h^{3,1}_{BC}(X) = h^{1,3}_{BC}(X) = h^{0,2}(X) , $$
$$ h^{2,2}_{BC}(X) = 2\, h^{2,1}_{BC}(X) - 2\, h^{0,1}(X) + 2 , $$
$$ h^{3,2}_{BC}(X) = h^{2,3}_{BC}(X) = h^{0,2}(X) + 1 + h^{2,0}(X) , $$
$$ h^{3,3}_{BC}(X) = 1 . $$
Here, we have denoted $h^{p,q}_{BC}(X):=\dim_{\C} H^{p,q}_{BC}(X)$, and we recall that $h^{p,q}_A(X):=\dim_{\C} H^{p,q}_{A}(X)=h^{3-q,3-p}
_{BC}(X)$ \cite[page 10]{schweitzer}.
At the end, he gets the dimensions of the Bott-Chern and Aeppli cohomologies in terms just of $h^{1,1}_{BC}(X)$, $h^{2,2}_{BC}(X)$, $h^{2,1}_{BC}(X)$, $h^{2,0}(X)$, $h^{0,2}(X)$.

Second, if we endow $S^6$ with a hypothetical complex structure and we construct the blowing-up at one point, then we get an exotic complex structure on a manifold diffeomorphic to $\mathbb{CP}^3$ (which is non-K\"ahler \cite{hirzebruch-kodaira}, not even Moishezon \cite{peternell}). Cohomological properties of hypothetical exotic complex structures on $\mathbb{CP}^3$ are investigated in \cite{brown}.

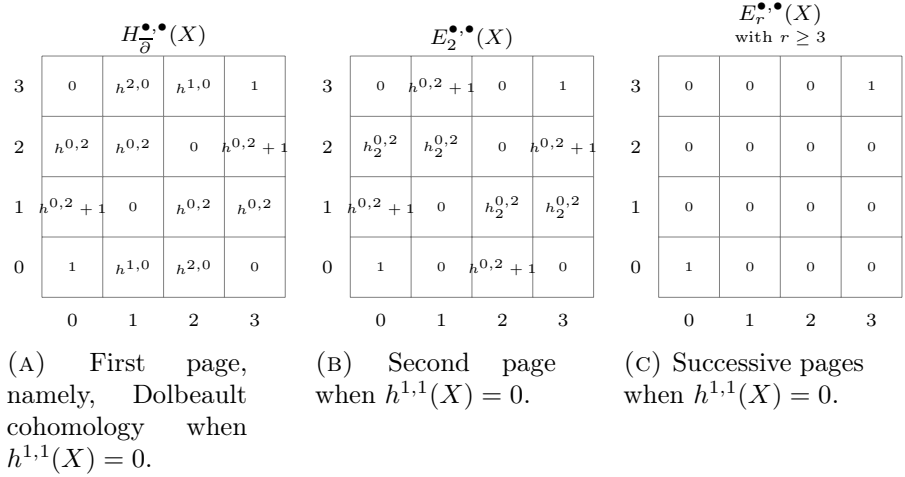
\begin{figure}
    \centering
    \begin{subfigure}[t]{0.25\textwidth}

\begin{tikzpicture}
\newcommand\un{.8}

\draw[help lines, step=\un] (0,0) grid (4*\un,4*\un);

\begingroup\makeatletter\def\f@size{8}\check@mathfonts
\node at (2*\un,4.3*\un) {$H^{\bullet,\bullet}_{\overline\partial}(X)$};
\endgroup

\begingroup\makeatletter\def\f@size{7}\check@mathfonts
\foreach \x in {0,...,3}
  \node at (\un*.5+\un*\x,-.3) {\x};
\foreach \y in {0,...,3}
  \node at (-.3,\un*.5+\un*\y) {\y};
\endgroup

\coordinate (A00) at (0*\un+1/2*\un, 0*\un+1/2*\un);
\coordinate (A01) at (0*\un+1/2*\un, 1*\un+1/2*\un);
\coordinate (A10) at (1*\un+1/2*\un, 0*\un+1/2*\un);
\coordinate (A20) at (2*\un+1/2*\un, 0*\un+1/2*\un);
\coordinate (A11) at (1*\un+1/2*\un, 1*\un+1/2*\un);
\coordinate (A02) at (0*\un+1/2*\un, 2*\un+1/2*\un);
\coordinate (A30) at (3*\un+1/2*\un, 0*\un+1/2*\un);
\coordinate (A21) at (2*\un+1/2*\un, 1*\un+1/2*\un);
\coordinate (A12) at (1*\un+1/2*\un, 2*\un+1/2*\un);
\coordinate (A03) at (0*\un+1/2*\un, 3*\un+1/2*\un);
\coordinate (A31) at (3*\un+1/2*\un, 1*\un+1/2*\un);
\coordinate (A22) at (2*\un+1/2*\un, 2*\un+1/2*\un);
\coordinate (A13) at (1*\un+1/2*\un, 3*\un+1/2*\un);
\coordinate (A32) at (3*\un+1/2*\un, 2*\un+1/2*\un);
\coordinate (A23) at (2*\un+1/2*\un, 3*\un+1/2*\un);
\coordinate (A33) at (3*\un+1/2*\un, 3*\un+1/2*\un);

\begingroup\makeatletter\def\f@size{4}\check@mathfonts
\node at (A00) {$1$};
\node at (A01) {$h^{0,2}+1$};
\node at (A10) {$h^{1,0}$};
\node at (A20) {$h^{2,0}$};
\node at (A11) {$0$};
\node at (A02) {$h^{0,2}$};
\node at (A30) {$0$};
\node at (A21) {$h^{0,2}$};
\node at (A12) {$h^{0,2}$};
\node at (A03) {$0$};
\node at (A31) {$h^{0,2}$};
\node at (A22) {$0$};
\node at (A13) {$h^{2,0}$};
\node at (A32) {$h^{0,2}+1$};
\node at (A23) {$h^{1,0}$};
\node at (A33) {$1$};
\endgroup

\end{tikzpicture}
\caption{First page, namely, Dolbeault cohomology when $h^{1,1}(X)=0$.}
\label{fig:dolbeault-h11zero}
    \end{subfigure}
    \qquad
    \begin{subfigure}[t]{0.25\textwidth}
\begin{tikzpicture}
\newcommand\un{.8}

\draw[help lines, step=\un] (0,0) grid (4*\un,4*\un);

\begingroup\makeatletter\def\f@size{8}\check@mathfonts
\node at (2*\un,4.3*\un) {$E^{\bullet,\bullet}_{2}(X)$};
\endgroup

\begingroup\makeatletter\def\f@size{7}\check@mathfonts
\foreach \x in {0,...,3}
  \node at (\un*.5+\un*\x,-.3) {\x};
\foreach \y in {0,...,3}
  \node at (-.3,\un*.5+\un*\y) {\y};
\endgroup

\coordinate (A00) at (0*\un+1/2*\un, 0*\un+1/2*\un);
\coordinate (A01) at (0*\un+1/2*\un, 1*\un+1/2*\un);
\coordinate (A10) at (1*\un+1/2*\un, 0*\un+1/2*\un);
\coordinate (A20) at (2*\un+1/2*\un, 0*\un+1/2*\un);
\coordinate (A11) at (1*\un+1/2*\un, 1*\un+1/2*\un);
\coordinate (A02) at (0*\un+1/2*\un, 2*\un+1/2*\un);
\coordinate (A30) at (3*\un+1/2*\un, 0*\un+1/2*\un);
\coordinate (A21) at (2*\un+1/2*\un, 1*\un+1/2*\un);
\coordinate (A12) at (1*\un+1/2*\un, 2*\un+1/2*\un);
\coordinate (A03) at (0*\un+1/2*\un, 3*\un+1/2*\un);
\coordinate (A31) at (3*\un+1/2*\un, 1*\un+1/2*\un);
\coordinate (A22) at (2*\un+1/2*\un, 2*\un+1/2*\un);
\coordinate (A13) at (1*\un+1/2*\un, 3*\un+1/2*\un);
\coordinate (A32) at (3*\un+1/2*\un, 2*\un+1/2*\un);
\coordinate (A23) at (2*\un+1/2*\un, 3*\un+1/2*\un);
\coordinate (A33) at (3*\un+1/2*\un, 3*\un+1/2*\un);

\begingroup\makeatletter\def\f@size{4}\check@mathfonts
\node at (A00) {$1$};
\node at (A01) {$h^{0,2}+1$};
\node at (A10) {$0$};
\node at (A20) {$h^{0,2}+1$};
\node at (A11) {$0$};
\node at (A02) {$h^{0,2}_2$};
\node at (A30) {$0$};
\node at (A21) {$h^{0,2}_2$};
\node at (A12) {$h^{0,2}_2$};
\node at (A03) {$0$};
\node at (A31) {$h^{0,2}_2$};
\node at (A22) {$0$};
\node at (A13) {$h^{0,2}+1$};
\node at (A32) {$h^{0,2}+1$};
\node at (A23) {$0$};
\node at (A33) {$1$};
\endgroup

\end{tikzpicture}
\caption{Second page when $h^{1,1}(X)=0$.}
\label{fig:E2-h11zero}
    \end{subfigure}
    \qquad
    \begin{subfigure}[t]{0.25\textwidth}
\begin{tikzpicture}
\newcommand\un{.8}

\draw[help lines, step=\un] (0,0) grid (4*\un,4*\un);

\begingroup\makeatletter\def\f@size{8}\check@mathfonts
\node at (2*\un,4.7*\un) {$E^{\bullet,\bullet}_{r}(X)$};
\endgroup
\begingroup\makeatletter\def\f@size{6}\check@mathfonts
\node at (2*\un,4.3*\un) {with $r\geq 3$};
\endgroup

\begingroup\makeatletter\def\f@size{7}\check@mathfonts
\foreach \x in {0,...,3}
  \node at (\un*.5+\un*\x,-.3) {\x};
\foreach \y in {0,...,3}
  \node at (-.3,\un*.5+\un*\y) {\y};
\endgroup

\coordinate (A00) at (0*\un+1/2*\un, 0*\un+1/2*\un);
\coordinate (A01) at (0*\un+1/2*\un, 1*\un+1/2*\un);
\coordinate (A10) at (1*\un+1/2*\un, 0*\un+1/2*\un);
\coordinate (A20) at (2*\un+1/2*\un, 0*\un+1/2*\un);
\coordinate (A11) at (1*\un+1/2*\un, 1*\un+1/2*\un);
\coordinate (A02) at (0*\un+1/2*\un, 2*\un+1/2*\un);
\coordinate (A30) at (3*\un+1/2*\un, 0*\un+1/2*\un);
\coordinate (A21) at (2*\un+1/2*\un, 1*\un+1/2*\un);
\coordinate (A12) at (1*\un+1/2*\un, 2*\un+1/2*\un);
\coordinate (A03) at (0*\un+1/2*\un, 3*\un+1/2*\un);
\coordinate (A31) at (3*\un+1/2*\un, 1*\un+1/2*\un);
\coordinate (A22) at (2*\un+1/2*\un, 2*\un+1/2*\un);
\coordinate (A13) at (1*\un+1/2*\un, 3*\un+1/2*\un);
\coordinate (A32) at (3*\un+1/2*\un, 2*\un+1/2*\un);
\coordinate (A23) at (2*\un+1/2*\un, 3*\un+1/2*\un);
\coordinate (A33) at (3*\un+1/2*\un, 3*\un+1/2*\un);

\begingroup\makeatletter\def\f@size{4}\check@mathfonts
\node at (A00) {$1$};
\node at (A01) {$0$};
\node at (A10) {$0$};
\node at (A20) {$0$};
\node at (A11) {$0$};
\node at (A02) {$0$};
\node at (A30) {$0$};
\node at (A21) {$0$};
\node at (A12) {$0$};
\node at (A03) {$0$};
\node at (A31) {$0$};
\node at (A22) {$0$};
\node at (A13) {$0$};
\node at (A32) {$0$};
\node at (A23) {$0$};
\node at (A33) {$1$};
\endgroup

\end{tikzpicture}
\caption{Successive pages when $h^{1,1}(X)=0$.}
\label{fig:E3-h11zero}
    \end{subfigure}
\caption{Dimensions of the Fr\"olicher spectral sequence of a hypothetical complex structure on $S^6$ when $h^{1,1}(X)=0$. With respect to the bi-graduation $\wedge^{p,q}X$, the horizontal axis represents $p$ and the vertical axis represents $q$.
The $h^{0,2}$, $h^{1,0}$, $h^{2,0}$, $h^{0,2}_2$ are unknown non-negative integers satisfying $h^{2,0}=h^{1,0}+h^{0,2}+1$.}
    \label{fig:frolicher-h11zero}
\end{figure}

Third, Gabor Etesi \cite{etesi} (see Remark 5.1 (3) in the arXiv version 8, \texttt{arXiv:math/0505634v8}) expects that $h^{0,1}(X) = h^{1,1}(X) = 1$ and $h^{0,2}(X) =  h^{1,2}(X) =  h^{1,0}(X) = h^{2,0}(X) = 0$. Bott-Chern cohomology in this case is considered in \cite[Section 3.2]{mchugh-1}. Under these assumptions, the Fr\"olicher spectral sequence degenerates at the second page.
On the other hand, the case $h^{1,1}(X)=0$ is investigated in \cite{ugarte}. Indeed, by the structure of $E_2^{\bullet,\bullet}(X)$ under the assumption $h^{1,1}(X)=0$, he gets that $h^{1,2}(X) = h^{0,2}(X)$, $h_2^{0,1}(X) = h^{0,1}(X)$, and then also $h^{0,1}(X) = h^{0,2}(X) + 1$, $h^{2,0}(X) = h^{1,0}(X) + h^{0,2}(X) + 1$. At the end, the Hodge numbers are completely determined by the non-negative integers $h^{1,0}=h^{1,0}(X)\in\{0,1\}$ and $h^{0,2}=h^{0,2}(X)$, and the Fr\"olicher spectral sequence is determined by $h^{1,0}$, $h^{0,2}$, and the non-negative integer $h^{0,2}_2=h^{0,2}_2(X)$, see Figure \ref{fig:frolicher-h11zero}.
Note in particular that $E_1^{\bullet,\bullet}(X) \neq E_2^{\bullet,\bullet}(X) \neq E_3^{\bullet,\bullet}(X) = E_\infty^{\bullet,\bullet}(X)$ if $h^{1,0}(X)=1$.

\bigskip

\noindent{\sl Acknowledgments.}
This note has been written for the MAM1 ``(Non)-existence of complex structures on $S^6$'' held in Marburg on March 27th--30th, 2017, \url{http://www.mathematik.uni-marburg.de/~agricola/Hopf2017/}.
The author warmly thanks the organizers for the kind hospitality, and also all the participants for the environment they contributed to.
Thanks in particular to Giovanni Bazzoni for his help and support during the preparation of the talk and of the paper.
The author would like to thank S\"onke Rollenske, Luis Ugarte, Claire Voisin for helpful and interesting discussions, and Andy McHugh for the references \cite{mchugh-1, mchugh-2, brown} and for several useful discussions.

\section{Preliminaries}
We denote by $X$ the manifold $S^6$ endowed with a hypothetical complex structure.

\subsection{Dolbeault cohomology and Fr\"olicher spectral sequence}

We recall that the {\em Dolbeault cohomology} is
$$ H^{p,q}_{\overline\partial}(X) = H^q(X;\Omega_{X}^p) $$
where $\Omega_{X}^p$ denotes the sheaf of germs of holomorphic $p$-forms over $X$.
Denote the Hodge numbers by $h^{p,q}(X):=\dim_\C H^{p,q}_{\overline\partial}(X)$.
By elliptic Hodge theory \cite{hodge}, the Dolbeault cohomology groups are isomorphic as vector spaces to the kernel of a $2$nd order elliptic operator, therefore the Hodge numbers are finite. Moreover, {\em Serre duality} holds: as real vector spaces,
$$ H^{p,q}_{\overline\partial}(X) \simeq H^{3-p,3-q}_{\overline\partial}(X) . $$
An explicit isomorphism can be costructed as follows. Fix a Hermitian metric $g$ on $X$, and consider the induced volume form. Then we have a (positive-definite) Hermitian product on $H^{\bullet,\bullet}_{\overline\partial}(X)$: the {\em $L^2$-pairing} $\langle \_ \vert \_ \rangle$. It can be represented by the $\C$-linear {\em Hodge-star-operator} $*\colon \wedge^{p,q}X\to\wedge^{3-q,3-p}X$ as
$$ \langle \alpha \vert \beta \rangle = \int_X \alpha \wedge \overline{*\beta} . $$
Then $\overline*\colon H^{p,q}_{\overline\partial}(X) \stackrel{\simeq}{\to} H^{3-q,3-p}_{\overline\partial}X$ is a $\C$-anti-linear isomorphism.

The double complex $\left( \wedge^{\bullet,\bullet}X, \partial, \overline\partial \right)$ gives a natural filtration (see {\itshape e.g.} \cite[Section 2.4]{mccleary})
$$ F^p(\wedge^\bullet S^6 \tens{\R} \C ) := \bigoplus_{\substack{r+s=\bullet \\ r\geq p}} \wedge^{r,s} X . $$
The associated spectral sequence is called {\em Fr\"olicher spectral sequence} \cite{frolicher}:
$$ E_1^{p,q}(X) \simeq H^{p,q}_{\overline\partial}(X) \Rightarrow H^{p+q}(S^6;\C) . $$
Denote $h_r^{p,q}(X):=\dim_\C E^{p,q}_{r}(X)$ for $r\in\mathbb{N}\cup\{\infty\}$ with $r\geq1$. Note that, for any $r,p,q,k$,
$$ h_{r+1}^{p,q}(X) \leq h_r^{p,q}(X) \qquad \text{ and } \qquad \sum_{p+q=k} h^{p,q}_\infty(X) = b_k(S^6) . $$
In particular, $h_r^{p,q}$ are finite, for $r\geq1$.

\begin{figure}[h]
    \centering
    \begin{subfigure}[t]{0.3\textwidth}
\centering
\begin{tikzpicture}
\newcommand\un{1}
\newcommand{\raggio}{2*\un pt}

\begingroup\makeatletter\def\f@size{6}\check@mathfonts
\node at (0*\un,3*\un) {$\cdots$};
\draw[->] (0*\un+1/4*\un,3*\un) -- (1*\un,3*\un) ;
\fill (1*\un+1/4*\un,3*\un) circle (\raggio);
\draw[->] (1*\un+1/4*\un,2*\un+1/4*\un) -- (1*\un+1/4*\un,3*\un-1/4*\un) ;
\fill (1*\un+1/4*\un,2*\un) circle (\raggio);
\draw[->] (1*\un+2/4*\un,2*\un) -- (2*\un,2*\un) ;
\fill (2*\un+1/4*\un,2*\un) circle (\raggio);
\draw[->] (2*\un+1/4*\un,1*\un+1/4*\un) -- (2*\un+1/4*\un,2*\un-1/4*\un) ;
\node at (2*\un+1/4*\un,1*\un) {$\vdots$};
\endgroup
\end{tikzpicture}

        \caption{Zigzags.}
        \label{subfig:zigzag}
    \end{subfigure}
    \qquad
   \begin{subfigure}[t]{0.3\textwidth}
\centering
\begin{tikzpicture}
\newcommand\un{1}
\newcommand{\raggio}{2*\un pt}

\begingroup\makeatletter\def\f@size{6}\check@mathfonts
\fill (0*\un,1*\un) circle (\raggio);
\fill (0*\un,0*\un) circle (\raggio);
\fill (1*\un,0*\un) circle (\raggio);
\fill (1*\un,1*\un) circle (\raggio);
\draw[->] (0*\un+1/4*\un,1*\un) -- (1*\un-1/4*\un,1*\un) ;
\draw[->] (0*\un+1/4*\un,0*\un) -- (1*\un-1/4*\un,0*\un) ;
\draw[->] (0*\un,0*\un+1/4*\un) -- (0*\un,1*\un-1/4*\un) ;
\draw[->] (1*\un,0*\un+1/4*\un) -- (1*\un,1*\un-1/4*\un) ;
\endgroup
\end{tikzpicture}
        \caption{Squares.}
        \label{subfig:square}
    \end{subfigure}
    \caption{Zigzags and squares of double complexes.}
    \label{fig:zigzag-square}
\end{figure}
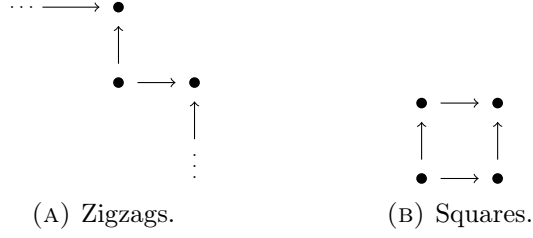

We give a heuristic picture, inspired by the MathOverflow discussion at \url{http://mathoverflow.net/questions/25723/} by Greg Kuperberg, see also \url{http://mathoverflow.net/questions/86947/}. See also \cite{angella-3}.
A double complex can be decomposed into direct sum of {\em zigzags} (Figure \ref{subfig:zigzag}) and {\em squares of isomorphisms} (Figure \ref{subfig:square}), see also \cite{deligne-griffiths-morgan-sullivan}. The {\em length} of a zigzag is defined as the number of non-zero objects in it. The squares do not contribute to cohomology. Even-length zigzags of length $2\ell$ do not contribute to the de Rham cohomology: they are killed at the $(\ell+1)$th page of the Fr\"olicher spectral sequence. Odd length zigzags gives contribution $1$ at every page of the spectral sequence, at one of the two ends of the zigzag itself. So, the Fr\"olicher spectral sequence degenerates at the first page if there is no zigzag of positive even length. In this case, the existence of zigzags of odd-length greater than one is then the obstruction to the Hodge decomposition (that is, to Dobeault cohomology providing a Hodge structure on the de Rham cohomology). On the other side, the maximum length of a even zigzag is $2\,\dim_\C X$, which means that the Fr\"olicher spectral sequence degenerate at most at the $(\dim_\C X+1)$th page. If there is no holomorphic top-degree form, then the maximum length of a zigzag is actually $2\,(\dim_\C X-1)$, which means that the Fr\"olicher spectral sequence degenerates at most at the $(\dim_\C X)$th page.

Therefore, by dimension reasons, the Fr\"olicher spectral sequence of $X$ converges at most at the $4$th page, being $4=\dim_\C X+1$. This means that, other than
$$ E_1^{p,q}\simeq H^{p,q}_{\overline\partial}(X) \qquad \text{ and } \qquad \bigoplus_{p+q=k}E_{\infty}^{p,q}\simeq H^k(S^6;\C) , $$
we have the additional finite-dimensional invariants:
\begin{eqnarray*}
E_2^{p,q} &\simeq& \frac{\ker (d_1\colon E_1^{p,q}\to E_1^{p+1,q})}{\im (d_1 \colon E_1^{p-1,q}\to E_1^{p,q})} \simeq H^p(H_{\overline\partial}^{\bullet,q},\partial), \\[5pt]
E_3^{p,q} &\simeq& \frac{\ker (d_2\colon E_2^{p,q}\to E_2^{p+2,q-1})}{\im (d_2 \colon E_2^{p-2,q+1}\to E_2^{p,q})} , \\[5pt]
E_4^{p,q} &\simeq& \frac{\ker (d_3\colon E_3^{p,q}\to E_3^{p+3,q-2})}{\im (d_3 \colon E_3^{p-3,q+2}\to E_3^{p,q})} .
\end{eqnarray*}
An explicit description of the above terms can be found in \cite[Theorem 1, Theorem 3]{cordero-fernandez-ugarte-gray}:
$$ E^{p,q}_r \simeq \frac{\mathcal{X}^{p,q}_r}{\mathcal{Y}^{p,q}_r} , $$
where, for $r=1$,
$$
 \mathcal{X}^{p,q}_1 := \left\{\alpha\in\wedge^{p,q}X \;:\; \overline\partial \alpha=0 \right\} ,
 \qquad
 \mathcal{Y}^{p,q}_1 := \overline\partial \wedge^{p,q-1}X ,
$$
and, for $r\geq 2$,
\begin{eqnarray*}
 \mathcal{X}^{p,q}_r &:=& \left\{\alpha^{p,q}\in\wedge^{p,q}X \;:\; \overline\partial \alpha^{p,q}=0 \text{ and, for any }i\in\{1, \ldots, r-1\}, \right. \\[5pt]
 && \left. \text{ there exists } \alpha^{p+i,q-i}\in\wedge^{p+i,q-i}X \right.\\[5pt]
 && \left. \text{ such that }\partial\alpha^{p+i-1,q-i+1}+\overline\partial\alpha^{p+i,q-i}=0 \right\} , \\[10pt]
 \mathcal{Y}^{p,q}_r &:=& \left\{\partial\beta^{p-1,q}+\overline\partial\beta^{p,q-1}\in\wedge^{p,q}X \;:\; \text{for any }i\in\{2, \ldots, r-1\}, \right. \\[5pt]
 && \left. \text{ there exists } \beta^{p-i,q+i-1}\in\wedge^{p-i,q+i-1}X \right. \\[5pt]
 && \left. \text{ such that } \partial\beta^{p-i,q+i-1}+\overline\partial\beta^{p-i+1,q+i-2}=0 \text{ and }\overline\partial\beta^{p-r+1, q+r-2}=0\right\} ,
\end{eqnarray*}
and, for any $r\geq 1$, the map $d_r\colon E^{p,q}_r\to E^{p+r, q-r+1}_r$ is given by
$$ d_r [\alpha^{p,q}] := [\partial\alpha^{p+r-1, q-r+1}] . $$

\section{Dolbeault cohomology}
In this section, we derive the Hodge numbers of a hypothetical complex structure on $S^6$.
By Serre duality, we have
$$ h^{3,3}(X) = h^{0,0}(X), \qquad h^{0,3}(X) = h^{3,0}(X), $$
$$ h^{3,1}(X) = h^{0,2}(X), \qquad h^{3,2}(X) = h^{0,1}(X), $$
$$ h^{1,3}(X) = h^{2,0}(X), \qquad h^{2,1}(X) = h^{1,2}(X), $$
$$ h^{2,2}(X) = h^{1,1}(X), \qquad h^{2,3}(X) = h^{1,0}(X). $$
We compute $h^{0,0}(X)=1$, $h^{3,0}(X)=0$, we express $h^{0,1}(X)$ in terms of the unknown non-negative integer $h^{0,2}(X)$ getting $h^{0,1}(X)\geq 1$, and we derive a relation between $h^{2,0}(X)$, $h^{1,1}(X)$, $h^{1,0}(X)$, $h^{1,2}(X)$ getting that either $h^{1,1}(X)$ or $h^{2,0}(X)$ is non-zero.
As for holomorphic $1$-forms, we can further prove that $h^{1,0}(X)\leq h^{2,0}(X)$. In Section \ref{subsec:further-2}, we will also derive the inequality $h^{1,1}(X) \geq h^{1,2}(X) - h^{0,2}(X)$.

We summarize the informations from the equations \eqref{eq:h00}, \eqref{eq:h30}, \eqref{eq:h01-h02}, \eqref{eq:h20-h11-h10-h12}, \eqref{eq:h10-h20}, \eqref{eq:h10-1}, \eqref{eq:h11-ugarte}, and Serre duality, in the Hodge diamond in Figure \ref{fig:dolbeault}.

\subsection{No non-zero holomorphic $3$-forms}
First of all, we notice that, $X$ being connected, then
\begin{equation}\label{eq:h00}
h^{0,0}(X) = 1 .
\end{equation}

We compute now the other corners of the Hodge diamond, that is, $h^{3,0}(X)$ and $h^{0,3}(X)$; note that they are equal by Serre duality.
We claim that $X$ does not have any non-zero holomorphic $3$-forms \cite[Lemma 4]{gray}, that is, $H^{3,0}_{\overline\partial}(X)=H^0(X;\Omega_X^3)=0$.
More precisely, we prove that there is a natural injective map $H^{3,0}_{\overline\partial}(X)\to H^3(S^6;\C)$, (more in general, there is a natural injective map $H^{n,0}_{\overline\partial}(X)\to H^n(X;\C)$ on a compact complex manifold of complex dimension $n$,) see {\itshape e.g.} \cite[Lemma 2.2]{brown}, from which we get the claim by using $H^3(S^6;\C)=0$. Let $[\beta]\in H^{3,0}_{\overline\partial}(X)$, that is, $\beta\in\wedge^{3,0}X$ such that $\overline\partial\beta=0$. By dimension reasons, also $d\beta=0$. Then $[\beta]\in H^3(S^6;\C)$. The map $H^{3,0}_{\overline\partial}(X) \ni [\beta] \mapsto [\beta] \in H^3(S^6;\C)$ is well-defined, since $\im \overline\partial \cap \wedge^{3,0}X = 0$.
Assume now that $[\beta]\mapsto 0$, that is, $\beta=d\eta$ for some $\eta\in\wedge^2S^6\tens{\R}\C$. Fix any Hermitian metric on $X$. Note that the associated ($\C$-linear) Hodge-star-operator $*\colon \wedge^{p,q}X\to\wedge^{3-q,3-p}$ is $\mathcal{C}^\infty(S^6;\C)$-linear. Note also that, with respect to a unitary local co-frame $\{\varphi^1,\varphi^2,\varphi^3\}$ for $\wedge^{1,0}X$, the form $\varphi^1\wedge\varphi^2\wedge\varphi^3$ is an eigenvector for $*$. This implies that $\overline{*\beta}=\overline\beta$.
By applying Stokes theorem, if $\beta=d\eta$ is exact, we get
\begin{eqnarray*}
\|\beta\|^2 &=& \langle \beta \vert \beta \rangle = \int_X \beta \wedge \overline{*\beta} = \int_X \beta \wedge \overline\beta \\[5pt]
&=& \int_X d\eta\wedge\overline{d\eta} = \int_X d(\eta\wedge\overline{d\eta}) = 0 .
\end{eqnarray*}
Since the $L^2$-pairing is non-degenerate, we get $\beta=0$, which proves the statement.
Summing up, we have computed
\begin{equation}\label{eq:h30}
h^{3,0}(X) = 0 .
\end{equation}

\subsection{The "surprising fact" that $h^{0,1}(X)$ is non-zero}
We prove now that \cite[Theorem 5]{gray}
$$ h^{0,1}(X) \geq 1 . $$
This is reported as a "surprising fact" by A. Gray: although $S^6$ is simply-connected, there exists a non-zero $\overline\partial$-closed non-$\overline\partial$-exact $(0,1)$-form that could be hopefully "interpreted geometrically", \cite[page 251]{gray}.
We give two proofs of this result: the first one by A. Gray \cite{gray} uses Atiyah and Singer index theorem; the second one by Joel Fine at \url{http://mathoverflow.net/questions/62492/} and by A. McHugh \cite[Theorem 2]{mchugh-1} uses the exponential sequence. Both of them rely on the property that $H^1(S^6;\Z)=H^2(S^6;\Z)=0$.

\begin{proof}[Proof (using Atiyah and Singer index theorem)]
Recall that the {\em Euler characteristic} is the index of the elliptic complex $(\wedge^\bullet S^6, d)$:
$$ \chi(S^6) := \sum_{k} (-1)^k b_k(S^6) = 1-0+0-0+0-0+1 = 2 , $$
where $b_k(S^6):=\dim_\R H^k(S^6;\R)$ denote the Betti numbers.
Note that the alternate sum does not change moving from $(E_r,d_r)$ to $(E_{r+1}=H(E_r,d_r), \, d_{r+1})$, that is, \cite[Theorem 1]{frolicher}
\begin{eqnarray*}
\chi &:=& \sum_k (-1)^k b_k(S^6) \\[5pt]
&=& \sum_{p,q} (-1)^{p+q} \dim_\C E^{p,q}_{r}(X) \\[5pt]
&=& \sum_{p,q} (-1)^{p+q} h^{p,q}(X)
\end{eqnarray*}
for any $r\in\{1,2,3,4\}$.

By taking into account the Serre duality, we can write:
\begin{eqnarray*}
\chi &=& \sum_{p,q} (-1)^{p+q} h^{p,q}(X) \\[5pt]
&=& 2 \, \underbrace{\sum_q (-1)^{q} h^{0,q}(X)}_{=:\chi_0(X)} + 2\, \sum_q (-1)^{1+q} h^{1,q} .
\end{eqnarray*}
Here,
$$ \chi_0(X) := \sum_q (-1)^{q} h^{0,q}(X) $$
is the {\em arithmetic genus} \cite{hirzebruch}, that is, the {\em analytic index} of the elliptic complex $(\wedge^{0,\bullet}X, \overline\partial)$.
By the Hirzebruch-Riemann-Roch theorem \cite{hirzebruch}, which is a special case of the Atiyah and Singer index theorem \cite[Theorem 2.12]{atiyah-singer}, the analytic index of the elliptic complex above is equal to the {\em topological index} of the elliptic operator $\overline\partial$, namely,
$$ \mathrm{ind}_t(\overline\partial) := \{\mathrm{ch}\,\sigma(\overline\partial)\,\mathrm{Td}(X)\}[TX] , $$
where $\mathrm{ch}\,\sigma(\overline\partial)$ denotes the Chern character of the symbol of the differential operator $\overline\partial$ and $\mathrm{Td}(X)$ is the Todd class of $X$ and $[TX]$ is the fundamental class of the tangent bundle. As a combination of Chern classes $c_1\in H^2(S^6;\Z)=0$ and $c_2\in H^4(S^6;\Z)=0$, we have
$$\mathrm{ind}_t(\overline\partial)=\mathrm{Td}(X)[X]=\frac{1}{24}c_1c_2[S^6]=0. $$
At the end, we get $\chi_0(X)=0$, that is,
$$ h^{0,0}(X)-h^{0,1}(X)+h^{0,2}(X)-h^{0,3}(X)=0 . $$
Taking into account \eqref{eq:h00}, and \eqref{eq:h30} and Serre duality, we get
\begin{equation}\label{eq:h01-h02}
h^{0,1}(X) = h^{0,2} + 1 \geq 1 ,
\end{equation}
which allow to write $h^{0,1}(X)$, (and $h^{3,2}(X)$, $h^{3,1}(X)$, by Serre duality,) in terms of the unknown non-negative integer $h^{0,2}=h^{0,2}(X)$.
\end{proof}

\begin{proof}[Proof (using the exponential sequence)]
Note that the property $h^{3,0}(X)=0$ assures that the canonical bundle $K_X$ is holomorphically non-trivial. Consider the exponential sequence
$$ 0 \to \Z \to \mathcal{O}_X \to \mathcal{O}_X^\times \to 0 . $$
The long exact sequence in cohomology gives, in particular,
$$ H^1(X;\Z) \to H^1(X;\mathcal{O}_X) \to H^1(X;\mathcal{O}_X^\times) \to H^2(X;\Z) . $$
Since $H^1(X;\Z)=0$ and $H^2(X;\Z)=0$, we get
$$ H^1(X;\mathcal{O}_X) \simeq H^1(X;\mathcal{O}_X^\times) , $$
and in particular $h^{0,1}(X)\neq0$, since $K_X \in H^1(X;\mathcal{O}_X^\times)$ is non-trivial.
\end{proof}

\subsection{The non-vanishing of Dolbeault cohomology at second degree}
Finally, by $2=\chi=2\, \sum_q (-1)^{1+q} h^{1,q}$, we get \cite[Proposition 3.1]{ugarte}
$$ h^{1,3}(X) + h^{1,1}(X) = h^{1,0}(X) + h^{1,2}(X) + 1 , $$
that is, by Serre duality,
\begin{equation}\label{eq:h20-h11-h10-h12}
h^{2,0}(X) + h^{1,1}(X) = h^{1,0}(X) + h^{1,2}(X) + 1 \geq 1.
\end{equation}
In particular, either $h^{2,0}(X)$ or $h^{1,1}(X)$ is non-zero.

\subsection{Further inequalities on the Hodge numbers -- I}
We derive here some further inequalities concerning the Hodge number $h^{1,0}(X)$, that is, the dimension of holomorphic $1$-forms, following \cite{huckleberry-kebekus-peternell, mchugh-1}.

We first derive the inequality \cite[Proposition 10.3]{huckleberry-kebekus-peternell}, \cite[Lemma 1]{mchugh-1}
\begin{equation}\label{eq:h10-h20}
h^{1,0}(X) \leq h^{2,0}(X) .
\end{equation}
We recall the proof by \cite{huckleberry-kebekus-peternell, mchugh-1}. Another explanation will be given in Section \ref{subsec:further-2} as in \cite[Remark 3.4]{ugarte}.
More precisely, we prove that the map $\partial \colon H^{1,0}_{\overline\partial}(X) \to H^{2,0}_{\overline\partial}(X)$ is injective, where $\partial$ denotes the induced map in Dolbeault cohomology. Indeed, let $\varphi$ be a $\overline\partial$-closed $(1,0)$-form such that $\partial[\varphi]=0$, that is, $[\partial\varphi]=0$. By degree reasons, it holds $\partial\varphi=0$. Therefore $d\varphi=0$. Since $b_1(S^6)=0$, then there exists a smooth function $f$ such that $\varphi=df$. By degree reasons, $\varphi=\partial f$ and $\overline\partial f=0$. This implies that $f$ is constant, and then $\varphi=0$.

\begin{remark}
Let us assume that the algebraic dimension $a(X):=\mathrm{tr}\,\mathrm{deg}_\C\mathcal{M}(X)$ of a hypothetical complex structure on $S^6$ is zero, as claimed in \cite[Corollary 2.2]{campana-demailly-peternell} and its Corrigendum; see the discussion in \cite{rollenske-mam1}.
In this case, we now prove the following bound for $h^{1,0}(X)$. It first appears in \cite[Proposition 10.3]{huckleberry-kebekus-peternell}, where it is attributed to Matei Toma (see also \cite[Lemma 2]{mchugh-1}).
We also thank S\"onke Rollenske for several remarks about it (see also \cite{rollenske-mam1}).
Under the assumption $a(X)=0$, it holds that:
\begin{equation}\label{eq:h10-1}
h^{1,0}(X) \leq 1 .
\end{equation}
Indeed, we first notice that, by degree reasons, Dolbeault cohomology $(\bullet,0)$-classes are identified with holomorphic $\bullet$-forms, which have a structure of algebra.
If we had $h^{1,0}(X) \geq 2$, then also $h^{2,0}(X)\geq h^{1,0}(X) \geq 2$ by \eqref{eq:h10-h20}. Suppose that $\alpha_1$ and $\alpha_2$ are independent holomorphic $1$-forms. Since by assumption the algebraic dimension $a(X):=\mathrm{tr}\,\mathrm{deg}_\C\mathcal{M}(X)$ is zero, that is, there is no non-constant meromorphic function on $X$, then $\alpha_1\wedge\alpha_2\neq0$. Then $\alpha_1\wedge\alpha_2$ is a holomorphic $2$-form; and there is another holomorphic $2$-form $\beta_2$ being independent with $\alpha_1\wedge\alpha_2$. Since $\beta_2$ is independent with $\alpha_1\wedge\alpha_2$ and using again that $a(X)=0$, then either $\alpha_1\wedge\beta_2$ or $\alpha_2\wedge\beta_2$ is a non-zero holomorphic $3$-form. This is impossible by \cite[Lemma 4]{gray}.
\end{remark}


\section{Fr\"olicher spectral sequence}
In this section, we compute the successive pages of the Fr\"olicher spectral sequence. We prove that it actually degenerates at the third page, as a consequence of the vanishing of holomorphic top-degree forms. Then we study the second page of the Fr\"olicher spectral sequence.
We summarize the informations in \eqref{eq:h203-h230}, \eqref{eq:h200-h233}, \eqref{eq:h210-h223}, \eqref{eq:h222-h211}, \eqref{eq:h2ug2}, \eqref{eq:h2ug3}, \eqref{eq:h2var} in the diamond in Figure \ref{fig:E2}.
In particular, we notice that a symmetry {\itshape à la} Serre holds also for $E_2^{\bullet,\bullet}$ \cite[page 174]{ugarte}. Moreover, either $h^{1,1}(X)\neq0$, or $E_2^{\bullet,\bullet}(X) \neq E_3^{\bullet,\bullet}(X) = E_\infty^{\bullet,\bullet}(X)$ where the Hodge numbers are completely determined by just $h^{1,0}(X)$ and $h^{0,2}(X)$, see Figure \ref{fig:frolicher-h11zero}.

\subsection{Degeneration at the third page}
We first prove that the Fr\"olicher spectral sequence of $X$ degenerates at the $3$rd page \cite[Lemma 2.1]{ugarte}.
This will follow by the property that $X$ has no holomorphic $3$-forms.

We recall that
$$ E^{p,q}_{r+1}(X) = \frac{\ker(d_{r} \colon E_{r}^{p,q}(X) \to E_{r}^{p+r,q-r+1}(X))}{\im (d_r\colon E_r^{p-r,q+r-1}(X) \to E_r^{p,q}(X))} . $$
Therefore, $E^{p,q}_{r+1}(X)$ are quotients of subgroups of $E_{r}^{p,q}(X)$. Since $\wedge^{p,q}X$ is non-zero only for $p,q\in\{0,1,2,3\}$, then $E^{p,q}_{r}(X)$ is non-zero only for $(p,q)\in\{0,1,2,3\}^2$, for any $r\geq1$.
In particular, $d_r=0$ for any $r\geq 4$, being $4=\dim_\C X+1$.
This means that $E_4^{\bullet,\bullet}(X)=E_\infty^{\bullet,\bullet}(X)$, and the Fr\"olicher spectral sequence degenerates at most at the $4$th page.
(More in general, on a compact complex manifold of complex dimension $n$, the same argument yields that the Fr\"olicher spectral sequence degenerates at most at the $(n+1)$th page.)

In this case, the only non-trivial maps that may appear in $( E_3^{\bullet,\bullet}(X), d_3)$ are:
$$ E_3^{0,2}(X) \to E_3^{3,0}(X) , \qquad E_3^{0,3}(X) \to E_3^{3,1}(X) . $$
Moreover $H^{3,0}_{\overline\partial}(X)=0$ by \eqref{eq:h30}, which implies also $H^{0,3}_{\overline\partial}(X)=0$ by Serre duality. Therefore $E^{3,0}_r(X) = E^{0,3}_r(X) = 0$ for any $r\geq1$. Then $d_3=0$, that is, $E_3^{\bullet,\bullet}(X) = E_{4}^{\bullet,\bullet}(X) = E_{\infty}^{\bullet,\bullet}(X)$, and the Fr\"olicher spectral sequence degenerates at most at the $3$th page.
(More in general, on a compact complex manifold of complex dimension $n$ with no non-zero holomorphic $n$-forms, the same argument yields that the Fr\"olicher spectral sequence degenerates at most at the $n$th page.)

Moreover, since $H^k(S^6;\R)=\bigoplus_{p+q=k}E_{\infty}^{p,q}=\bigoplus_{p+q=k}E_{3}^{p,q}$, then the only non-trivial dimensions, for $r\geq3$, are
\begin{equation}\label{eq:hr00-hr33}
h_r^{0,0}(X)=1 , \qquad h_r^{3,3}(X)=1,
\end{equation}
see Figure \ref{fig:E3}.

\subsection{Second page}

We compute now the $2$nd page of the Fr\"olicher spectral sequence \cite[Theorem 3.2]{ugarte}.

We first compute the four corners of the diamond.
Since, by \eqref{eq:h30} and Serre duality, we have $0=h^{3,0}(X)=h_1^{3,0}(X) \geq h_2^{3,0}(X)$ and $0 = h^{0,3}(X) = h_1^{0,3}(X) \geq h_2^{0,3}(X)$, then
\begin{equation}\label{eq:h203-h230}
h_2^{3,0}(X) = 0, \qquad h_2^{0,3}(X) = 0 .
\end{equation}
Since, by \eqref{eq:hr00-hr33} and \eqref{eq:h00} and Serre duality, we have $1 \leq h_3^{0,0}(X) \leq h_2^{0,0}(X) \leq h_1^{0,0}(X) = h^{0,0}(X) = 1$, where the last equality follows by the Maximum Principle, $X$ being compact, and $1 \leq h_3^{3,3}(X) \leq h_2^{3,3}(X) \leq h_1^{3,3}(X) = 1$, then
\begin{equation}\label{eq:h200-h233}
h_2^{0,0}(X) = 1, \qquad h_2^{3,3}(X) = 1 .
\end{equation}

The non-trivial maps in the complex $(E^{\bullet,\bullet}_2(X) , d_2\colon E^{p,q}_2(X) \to E^{p+2,q-1}_2(X) )$, by dimension reasons, are
$$ 0 \to E_2^{0,1}(X) \to E_2^{2,0}(X) \to 0 , \qquad 0 \to E_2^{1,1}(X) \to E_2^{3,0}(X) \to 0 , $$
$$ 0 \to E_2^{0,2}(X) \to E_2^{2,1}(X) \to 0 , \qquad 0 \to E_2^{1,2}(X) \to E_2^{3,1}(X) \to 0 , $$
$$ 0 \to E_2^{0,3}(X) \to E_2^{2,2}(X) \to 0 , \qquad 0 \to E_2^{1,3}(X) \to E_2^{3,2}(X) \to 0 , $$
while
$$ 0 \to E_2^{1,0}(X) \to 0 , \qquad 0 \to E_2^{2,3}(X) \to 0 . $$
Actually, by \eqref{eq:h203-h230}, we have
$$ 0 \to E_2^{1,1}(X) \to 0 , \qquad 0 \to E_2^{2,2}(X) \to 0 . $$
Since the cohomology of these maps is given by $E^{p,q}_3(X)$ that vanishes for any $(p,q)\not\in\{(0,0),(3,3)\}$, we get
\begin{equation}\label{eq:h210-h223}
h_2^{1,0}(X) = 0 , \qquad h_2^{2,3}(X) = 0 ,
\end{equation}
and
\begin{equation}\label{eq:h222-h211}
h_2^{1,1}(X) = 0 , \qquad h_2^{2,2}(X) = 0 ,
\end{equation}
and isomorphisms giving
\begin{equation}\label{eq:h2ug1a}
h_2^{0,1}(X) = h_2^{2,0}(X) ,
\qquad
h_2^{1,3}(X) = h_2^{3,2}(X) ,
\end{equation}
\begin{equation}\label{eq:h2ug1b}
h_2^{1,2}(X) = h_2^{3,1}(X) ,
\qquad
h_2^{0,2}(X) = h_2^{2,1}(X) .
\end{equation}

A further relation can be obtained by the following observation.
For $q\in\{0,1,2,3\}$, the complex
$$ 0 \to E_1^{0,q}(X) \stackrel{d_1}{\to} E_1^{1,q}(X) \stackrel{d_1}{\to} E_1^{2,q}(X) \stackrel{d_1}{\to} E_1^{3,q}(X) \to 0 . $$
has the same index as
$$ 0 \to E_2^{0,q}(X) \to E_2^{1,q}(X) \to E_2^{2,q}(X) \to E_2^{3,q}(X) \to 0 . $$
Then we get
\begin{eqnarray*}
\lefteqn{ h_1^{0,q}(X) - h_1^{1,q}(X) + h_1^{2,q}(X) - h_1^{3,q}(X) } \\[5pt]
&=& h_2^{0,q}(X) - h_2^{1,q}(X) + h_2^{2,q}(X) - h_2^{3,q}(X) .
\end{eqnarray*}
Varying $q\in\{0,1,2,3\}$ and using the known relations on the Hodge numbers, we have:
\begin{eqnarray*}
1-h^{1,0}(X)+h^{1,3}(X)
&=&
1+h_2^{2,0}(X)
, \\[5pt]
h^{0,2}(X)+1-h^{1,1}(X)+h^{1,2}(X)-h^{0,2}(X)
&=&
h_2^{0,1}(X)+h_2^{2,1}(X)-h_2^{3,1}(X)
, \\[5pt]
h^{0,2}(X)-h^{1,2}(X)+h^{1,1}(X)-h^{0,2}(X)-1
&=&
h_2^{0,2}(X)-h_2^{1,2}(X)-h_2^{3,2}(X)
, \\[5pt]
-h^{1,3}(X)+h^{1,0}(X)-1
&=&
-h_2^{1,3}(X)-1 .
\end{eqnarray*}
In particular, the first and the fourth equations give $h_2^{2,0}(X) = h_2^{1,3}(X)$, whence, by \eqref{eq:h2ug1a},
\begin{equation}\label{eq:h2ug2}
h_2^{0,1}(X) = h_2^{2,0}(X) = h_2^{1,3}(X) = h_2^{3,2}(X) .
\end{equation}
The second and third equations now give $h_2^{0,2}(X) = h_2^{1,2}(X)$, whence, by \eqref{eq:h2ug1b},
\begin{equation}\label{eq:h2ug3}
h_2^{2,1}(X) = h_2^{0,2}(X) = h_2^{1,2}(X) = h_2^{3,1}(X) .
\end{equation}
In particular, we notice that a duality {\itshape à la} Serre holds also for $E_2^{\bullet,\bullet}$, namely, for any $p,q$,
$$ h_2^{p,q} = h_2^{3-p,3-q} . $$

Again from the second equation, we have also the equality
\begin{equation}\label{eq:h2var}
h^{0,1}_2(X) = h^{1,2}(X) - h^{1,1}(X) + 1 \geq 0 .
\end{equation}
We notice that \cite[Corollary 3.3]{ugarte}: either $h^{1,1}(X)\neq 0$; or $h^{1,1}(X)=0$ then $h_2^{0,1}(X)\neq 0=h_3^{0,1}(X)$, that is, $E_2^{\bullet,\bullet}(X) \neq E_{\infty}^{\bullet,\bullet}(X)$.

\subsection{Further inequalities on the Hodge numbers -- II}\label{subsec:further-2}
By the vanishing of $E_2^{1,0}(X)$ and $E_2^{1,1}(X)$, we can deduce further inequalities on the Hodge numbers \cite[Remark 3.4]{ugarte}, \cite[Proposition 10.3]{huckleberry-kebekus-peternell}, \cite[Lemma 1]{mchugh-1}.

First, we give another evidence for the inequality \eqref{eq:h10-h20}.
Consider the complex
$$ \C = E_1^{0,0}(X) \stackrel{d_1=0}{\to} E_1^{1,0}(X) \stackrel{d_1}{\to} E_1^{2,0}(X) \stackrel{d_1}{\to} E_1^{3,0}(X) = 0 . $$
Since $E_2^{1,0}(X)=H^{1,0}(E_1^{\bullet,\bullet}(X),d_1)=0$, then the above complex is exact at $E_1^{1,0}(X)$, that is, the map $E_1^{1,0}(X) \to E_1^{2,0}(X)$ is injective. We get again \eqref{eq:h10-h20}:
$$ h^{1,0}(X) \leq h^{2,0}(X) . $$

Similarly, consider the complex
$$ 0 \to E_1^{0,1}(X) \stackrel{d_1^{0,1}}{\to} E_1^{1,1}(X) \stackrel{d_1^{1,1}}{\to} E_1^{2,1}(X) \stackrel{d_1^{2,1}}{\to} E_1^{3,1}(X) \to 0 . $$
By using that $E_2^{1,1}(X) = H^{1,1}(E_1^{\bullet,\bullet}(X), d_1) = 0$, we split it into the exact sequences
$$ 0 \to \ker d_1^{0,1} \to E_1^{0,1}(X) \stackrel{d_1^{0,1}}{\to} E_1^{1,1}(X) \stackrel{d_1^{1,1}}{\to} E_1^{2,1}(X) \to \coker d_1^{1,1} \to 0 , $$
$$ 0 \to E_2^{2,1}(X) = \frac{\ker d_1^{2,1}}{\im d_1^{1,1}} \to \coker d_1^{1,1} \stackrel{d_1^{2,1}}{\to} E_1^{3,1}(X) \to \coker d_1^{2,1} = E_2^{3,1}(X) \to 0 . $$
Then we get
\begin{eqnarray*}
h^{1,1}(X)
&=& \left( h^{0,1}(X) - \dim \ker d_1^{0,1} \right) + h^{2,1}(X) - \dim \coker d_1^{1,1} \\[5pt]
&\geq& h^{2,1}(X) - \left( h_2^{2,1}(X) + h^{3,1}(X) - h_2^{3,1}(X) \right) \\[5pt]
&=& h^{2,1}(X) - h^{3,1}(X) ,
\end{eqnarray*}
from which we finally get \cite[Remark 3.4]{ugarte}
\begin{equation}\label{eq:h11-ugarte}
h^{1,1}(X) \geq h^{1,2}(X) - h^{0,2}(X) .
\end{equation}

\section{Bott-Chern cohomology}
From the knowledge of the pages of the Fr\"olicher spectral sequence, we get the double complex of forms: its main structure is depicted it in Figure \ref{fig:S6}, see \cite{angella-3}.
Indeed the odd-length zigzags, which are the ones contributing to the de Rham cohomology, sit only at $E_\infty^{0,0}(X)$ and $E_\infty^{3,3}(X)$.

In particular, we can deduce informations also on the Bott-Chern and Aeppli cohomologies,
$$ H^{\bullet,\bullet}_{BC}(X) := \frac{\ker\partial\cap\ker\overline\partial}{\im\partial\overline\partial}, \qquad
H^{\bullet,\bullet}_{A}(X) := \frac{\ker\partial\overline\partial}{\im\partial+\im\overline\partial}. $$
Denote $h^{p,q}_{BC}(X):=\dim_\C H^{p,q}_{BC}(X)$ and $h^{p,q}_{A}(X):=\dim_\C H^{p,q}_{A}(X)$.
Recall that \cite[page 10]{schweitzer}, since $X$ is compact of complex dimension $3$, we have
$$ h^{p,q}_{A}(X) = h^{3-q,3-p}_{BC}(X) . $$ 

The Bott-Chern cohomology of $X$ is investigated in \cite{mchugh-1, mchugh-2} by Andrew McHugh. More precisely, he proves the following results \cite[Section 3]{mchugh-1}, \cite[Proposition 3.5]{mchugh-2}:
\begin{itemize}
\item $h^{0,0}_{BC}(X) = 1$: indeed \cite[Theorem 4]{mchugh-1}, $d$-closed functions on $X$ are locally-constant; since $X$ is connected, they are constant;
\item $h^{1,0}_{BC}(X) = 0$: indeed, $H^{1,0}_{BC}(X)=\wedge^{1,0}X \cap \ker\partial \cap \ker\overline\partial \simeq E_2^{1,0}(X)=0$. More concretely \cite[Lemma 7]{mchugh-1}, let $\varphi\in \wedge^{1,0}X$ such that $\partial\varphi=\overline\partial\varphi=0$; then also $d\varphi=0$ and $[\varphi]\in H^1_{dR}(X;\C)=0$; then there exists $f\in\mathcal{C}^{\infty}(X;\C)$ such that $\varphi=d f$; this means that $\varphi=\partial f$ and $\overline\partial f=0$; that is, $f$ is a holomorphic function on $X$ compact connected, whence constant by the Maximum Principle; we get that $\varphi=0$;
\item $h^{0,1}_{BC}(X) = 0$: indeed, conjugation induces an $\R$-linear isomorphism $H^{p,q}_{BC}(X)\simeq H^{q,p}_{BC}(X)$, for any bi-degree $(p,q)$;
\item $h^{1,1}_{BC}(X) = 2\,h^{0,1}(X)$: indeed \cite[Proposition 4.1]{mchugh-2}, by using that $b_1(S^6)=0$, we have that $H^{1,1}_{BC}(X) = \mathrm{im}\left( \overline\partial \colon H^{1,0}_{\partial}(X) \to H^{1,1}_{BC}(X) \right) + \mathrm{im}\left( \partial \colon H^{0,1}_{\overline\partial}(X) \to H^{1,1}_{BC}(X) \right)$; note now that $b_1(S^6)=0$ yields that this is in fact a direct sum; finally, the maps $\overline\partial \colon H^{1,0}_{\partial}(X) \to H^{1,1}_{BC}(X)$ and $\partial \colon H^{0,1}_{\overline\partial}(X) \to H^{1,1}_{BC}(X)$ are injective since $b_1(S^6)=0$;
\item $h^{2,0}_{BC}(X) = h^{2,0}(X)$: indeed \cite[Lemma 9]{mchugh-1}, we have an injective map $H^{2,0}_{BC}(X) = \wedge^{2,0}X \cap \ker\partial \cap \ker\overline\partial \hookrightarrow \wedge^{2,0}X \cap \ker\overline\partial = H^{2,0}_{\overline\partial}(X)$; on the other side, any $\overline\partial$-closed $(2,0)$-form $\varphi$ is also $\partial$-closed: this is because $\partial\varphi$ is a holomorphic $3$-form, but there is no non-zero holomorphic $3$-form \cite[Lemma 4]{gray}; then the natural map $H^{2,0}_{BC}(X) \to H^{2,0}_{\overline\partial}(X)$ induced by the identity is an isomorphism;
\item $h^{0,2}_{BC}(X) = h^{2,0}(X)$: by conjugation;
\item $h^{3,0}_{BC}(X) = 0$: we have \cite[Lemma 8]{mchugh-1} that $H^{3,0}_{BC}(X) = \wedge^{3,0}X \cap \ker\partial \cap \ker\overline\partial \hookrightarrow \wedge^{3,0}X \cap \ker\overline\partial = H^{3,0}_{\overline\partial}(X) = 0$;
\item $h^{0,3}_{BC}(X) = 0$: by conjugation;
\item $h^{1,3}_{BC}(X) = h^{0,2}(X)$: indeed \cite[Lemma 9]{mchugh-1}, consider the sequence
$$ 0 \to H^{0,2}_{\overline\partial}(X) \to H^{0,2}_{A}(X) \stackrel{\overline\partial}{\to} H^{0,3}_{BC}(X) ; $$
it is straightforward to check that it is exact, because of degree reasons; moreover, we know that $h^{0,3}_{BC}(X) = 0$; then we get that $h^{1,3}_{BC}(X)=h^{0,2}_{A}(X)=h^{0,2}_{\overline\partial}(X)$, where we have used Schweitzer duality;
\item $h^{2,2}_{BC}(X)=2\,h^{2,1}_{BC}(X)-2\,h^{0,1}(X)+2$: this follows by the computation of $h^{1,1}_{BC}(X)=2\,h^{0,1}(X)$ above and by the relation
$$ h^{1,1}_{BC}(X) + h^{2,2}_{BC}(X) = 2\, h^{2,1}_{BC}(X) + 2 , $$
that we prove now;
indeed \cite[page 13]{mchugh-1}, we have the exact sequence \cite[Theorem 3]{mchugh-1}
\begin{eqnarray*}
&& 0 \to H^{1,0}_{\overline\partial}(X) \to H^{1,0}_{A}(X) \\[5pt]
&& \stackrel{\overline\partial}{\to} H^{1,1}_{BC}(X) \to H^{1,1}_{\overline\partial}(X) \to H^{1,1}_{A}(X) \\[5pt]
&& \stackrel{\overline\partial}{\to} H^{1,2}_{BC}(X) \to H^{1,2}_{\overline\partial}(X) \to H^{1,2}_{A}(X) \\[5pt]
&& \stackrel{\overline\partial}{\to} H^{1,3}_{BC}(X) \to H^{1,3}_{\overline\partial}(X) \to H^{1,3}_{A}(X) \stackrel{\overline\partial}{\to} 0
\end{eqnarray*}
because of $b_1(S^6)=b_2(S^6)=b_3(S^6)=b_4(S^6)=0$. By using the previous computations, as well as the computations of the Hodge numbers, this gives
\begin{eqnarray*}
0 &=& h^{1,0}(X)-(h^{0,2}(X)+1+h^{2,0}(X))+h^{1,1}_{BC}(X) \\[5pt]
&& - h^{1,1}(X)+h^{2,2}_{BC}(X)-h^{2,1}_{BC}(X) \\[5pt]
&& +h^{1,2}(X)-h^{2,1}_{BC}(X)+h^{0,2}(X) \\[5pt]
&& -h^{2,0}(X)+h^{2,0}(X) .
\end{eqnarray*}
By using also \eqref{eq:h20-h11-h10-h12}, this completes the proof of the claim.
\item $h^{3,1}_{BC}(X) = h^{0,2}(X)$: by conjugation;
\item $h^{2,3}_{BC}(X) = h^{0,2}(X) + 1 + h^{2,0}(X)$: indeed \cite[Lemma 9]{mchugh-1}, consider the sequence
$$ 0 \to H^{0,1}_{\overline\partial}(X) \to H^{0,1}_{A}(X) \stackrel{\overline\partial}{\to} H^{0,2}_{BC}(X) \to H^{0,2}_{\overline\partial}(X) ; $$
it is straghtforward to check that it is exact, because of degree reasons; moreover, since $b_2(S^6)=0$, then the image of the map $H^{0,2}_{BC}(X) \to H^{0,2}_{\overline\partial}(X)$ is zero; then we get that $h^{2,3}_{BC}(X)=h^{0,1}_{A}(X)=h^{0,1}(X)+h^{0,2}_{BC}(X)=h^{0,2}(X)+1+h^{2,0}(X)$, where we have used the Schweitzer duality, \eqref{eq:h01-h02}, and the previous computation of $h^{2,0}_{BC}(X)=h^{2,0}(X)$;
\item $h^{3,2}_{BC}(X) = h^{0,2}(X) + 1 + h^{2,0}(X)$: by conjugation;
\item $h^{3,3}_{BC}(X) = 1$: by the Schweitzer duality, $H^{3,3}_{BC}(X) \simeq H^{0,0}_{A}(X) = \wedge^{0,0}X \cap \ker \partial\overline\partial$: this is the space of pluri-harmonic functions on $X$ compact connected, which are constant by the Maximum Principle; see also \cite[Theorem 4]{mchugh-1}.
\end{itemize}
Summing up, he can express the Bott-Chern cohomology dimensions in terms just of $h^{1,1}_{BC}(X)$, $h^{2,2}_{BC}(X)$, $h^{2,1}_{BC}(X)$, $h^{2,0}(X)$, $h^{0,2}(X)$, see Figure \ref{fig:BC-A}.

\begin{figure}
    \centering
    \begin{subfigure}[t]{0.42\textwidth}
\begin{tikzpicture}
\newcommand\un{1.5}

\draw[help lines, step=\un] (0,0) grid (4*\un,4*\un);

\begingroup\makeatletter\def\f@size{8}\check@mathfonts
\node at (2*\un,4.3*\un) {$H^{\bullet,\bullet}_{BC}(X)$};
\endgroup

\begingroup\makeatletter\def\f@size{7}\check@mathfonts
\foreach \x in {0,...,3}
  \node at (\un*.5+\un*\x,-.3) {\x};
\foreach \y in {0,...,3}
  \node at (-.3,\un*.5+\un*\y) {\y};
\endgroup

\coordinate (A00) at (0*\un+1/2*\un, 0*\un+1/2*\un);
\coordinate (A01) at (0*\un+1/2*\un, 1*\un+1/2*\un);
\coordinate (A10) at (1*\un+1/2*\un, 0*\un+1/2*\un);
\coordinate (A20) at (2*\un+1/2*\un, 0*\un+1/2*\un);
\coordinate (A11) at (1*\un+1/2*\un, 1*\un+1/2*\un);
\coordinate (A11a) at (1*\un+1/2*\un, 1*\un+1/2*\un+1/9*\un);
\coordinate (A11b) at (1*\un+1/2*\un, 1*\un+1/2*\un-1/9*\un);
\coordinate (A02) at (0*\un+1/2*\un, 2*\un+1/2*\un);
\coordinate (A30) at (3*\un+1/2*\un, 0*\un+1/2*\un);
\coordinate (A21) at (2*\un+1/2*\un, 1*\un+1/2*\un);
\coordinate (A12) at (1*\un+1/2*\un, 2*\un+1/2*\un);
\coordinate (A03) at (0*\un+1/2*\un, 3*\un+1/2*\un);
\coordinate (A31) at (3*\un+1/2*\un, 1*\un+1/2*\un);
\coordinate (A22) at (2*\un+1/2*\un, 2*\un+1/2*\un);
\coordinate (A22a) at (2*\un+1/2*\un, 2*\un+1/2*\un+1/9*\un);
\coordinate (A22b) at (2*\un+1/2*\un, 2*\un+1/2*\un-1/9*\un);
\coordinate (A13) at (1*\un+1/2*\un, 3*\un+1/2*\un);
\coordinate (A32) at (3*\un+1/2*\un, 2*\un+1/2*\un);
\coordinate (A23) at (2*\un+1/2*\un, 3*\un+1/2*\un);
\coordinate (A33) at (3*\un+1/2*\un, 3*\un+1/2*\un);

\begingroup\makeatletter\def\f@size{2}\check@mathfonts
\node at (A00) {$1$};
\node at (A01) {$0$};
\node at (A10) {$0$};
\node at (A20) {$h^{2,0}$};
\node at (A11) {$2h^{0,1}$};
\node at (A02) {$h^{2,0}$};
\node at (A30) {$0$};
\node at (A21) {$h^{2,1}_{BC}$};
\node at (A12) {$h^{2,1}_{BC}$};
\node at (A03) {$0$};
\node at (A31) {$h^{0,2}$};
\node at (A22a) {$2h^{2,1}_{BC}$};
\node at (A22b) {$- 2h^{0,1} + 2$};
\node at (A13) {$h^{0,2}$};
\node at (A32) {$h^{0,2} + 1 + h^{2,0}$};
\node at (A23) {$h^{0,2} + 1 + h^{2,0}$};
\node at (A33) {$1$};
\endgroup

\end{tikzpicture}
\caption{Bott-Chern cohomology.}
\label{fig:BC}
    \end{subfigure}
    \qquad\qquad
    \begin{subfigure}[t]{0.42\textwidth}
\begin{tikzpicture}
\newcommand\un{1.5}

\draw[help lines, step=\un] (0,0) grid (4*\un,4*\un);

\begingroup\makeatletter\def\f@size{8}\check@mathfonts
\node at (2*\un,4.3*\un) {$H^{\bullet,\bullet}_{A}(X)$};
\endgroup

\begingroup\makeatletter\def\f@size{7}\check@mathfonts
\foreach \x in {0,...,3}
  \node at (\un*.5+\un*\x,-.3) {\x};
\foreach \y in {0,...,3}
  \node at (-.3,\un*.5+\un*\y) {\y};
\endgroup

\coordinate (A00) at (0*\un+1/2*\un, 0*\un+1/2*\un);
\coordinate (A01) at (0*\un+1/2*\un, 1*\un+1/2*\un);
\coordinate (A10) at (1*\un+1/2*\un, 0*\un+1/2*\un);
\coordinate (A20) at (2*\un+1/2*\un, 0*\un+1/2*\un);
\coordinate (A11) at (1*\un+1/2*\un, 1*\un+1/2*\un);
\coordinate (A02) at (0*\un+1/2*\un, 2*\un+1/2*\un);
\coordinate (A30) at (3*\un+1/2*\un, 0*\un+1/2*\un);
\coordinate (A21) at (2*\un+1/2*\un, 1*\un+1/2*\un);
\coordinate (A12) at (1*\un+1/2*\un, 2*\un+1/2*\un);
\coordinate (A03) at (0*\un+1/2*\un, 3*\un+1/2*\un);
\coordinate (A31) at (3*\un+1/2*\un, 1*\un+1/2*\un);
\coordinate (A22) at (2*\un+1/2*\un, 2*\un+1/2*\un);
\coordinate (A13) at (1*\un+1/2*\un, 3*\un+1/2*\un);
\coordinate (A32) at (3*\un+1/2*\un, 2*\un+1/2*\un);
\coordinate (A23) at (2*\un+1/2*\un, 3*\un+1/2*\un);
\coordinate (A33) at (3*\un+1/2*\un, 3*\un+1/2*\un);

\begingroup\makeatletter\def\f@size{2}\check@mathfonts
\node at (A00) {$1$};
\node at (A01) {$h^{0,2} + 1 + h^{2,0}$};
\node at (A10) {$h^{0,2} + 1 + h^{2,0}$};
\node at (A20) {$h^{0,2}$};
\node at (A11a) {$2h^{2,1}_{BC}$};
\node at (A11b) {$- 2h^{0,1} + 2$};
\node at (A02) {$h^{0,2}$};
\node at (A30) {$0$};
\node at (A21) {$h^{2,1}_{BC}$};
\node at (A12) {$h^{2,1}_{BC}$};
\node at (A03) {$0$};
\node at (A31) {$h^{2,0}$};
\node at (A22) {$2h^{0,1}$};
\node at (A13) {$h^{2,0}$};
\node at (A32) {$0$};
\node at (A23) {$0$};
\node at (A33) {$1$};
\endgroup

\end{tikzpicture}
\caption{Aeppli cohomology.}
\label{fig:A}
    \end{subfigure}
\caption{Dimensions of the Bott-Chern and Aeppli cohomologies of a hypothetical complex structure on $S^6$, see \cite{mchugh-1, mchugh-2}. With respect to the bi-graduation $\wedge^{p,q}X$, the horizontal axis represents $p$ and the vertical axis represents $q$.
The $h^{2,1}_{BC}$, $h^{2,0}$, $h^{0,2}$ are unknown non-negative integers.}
    \label{fig:BC-A}
\end{figure}

\section{Open problems}

\begin{problem}[Sullivan-Barge Theorem for complex manifolds]
The Sullivan-Barge Theorem \cite{barge, sullivan} states that the necessary conditions for the realization of a Sullivan model by a compact simply-connected manifold are also sufficient.
Determine necessary and sufficient conditions for the realizability of a double complex as the double complex of a compact complex manifold.
\end{problem}

\begin{problem}[twisted Dolbeault cohomology of hypothetical complex manifold $X$ with underlying $S^6$]
Let $[\vartheta^{0,1}]$ be a non-trivial class in $H^{0,1}_{\overline\partial}(X)$. Study the cohomology of the twisted complex $(\wedge^{p,\bullet}X, \overline\partial-\vartheta^{0,1}\wedge\_)$
\end{problem}

\end{document}